%% file: SRMsurface.tex
\documentclass[a4paper,11pt,leqno]{article}
\usepackage{latexsym,amssymb,amsmath}
\usepackage{amsfonts}
\usepackage{eufrak}
\usepackage{graphicx}

\include{newcom1}

\newcommand{\la}{\langle}
\newcommand{\ra}{\rangle}
\newtheorem{example}[theorem]{Example}
\begin{document}

\title{ Minimal surfaces in contact\\ Sub - Riemannian manifolds.}
\author{Nataliya Shcherbakova}
\date{}
\maketitle

\begin{abstract}

In the present paper we consider generic Sub-Riemannian structures on the co-rank 1 non-holonomic vector 
distributions and introduce the associated canonical volume and ''horizontal'' area forms.
As in the classical case, the Sub-Riemannian minimal surfaces can be defined as the critical points of the 
'`horizontal''  area functional. We derive an intrinsic equation for minimal surfaces 
associated to a generic Sub-Riemannian structure of co-rank 1 in terms of the canonical volume form 
and the ``horizontal'' normal. The presented construction permits to describe the Sub-Riemannian minimal surfaces in a 
generic Sub-Riemannian manifold and can be easily generalized to the case of non-holonomic 
vector distributions of greater co-rank.

The case of contact vector distributions, in particular the $(2,3)$-case,  is studied more in detail.
In the latter case the geometry of the Sub-Riemannian minimal surfaces  is determined by the structure of their 
characteristic points (i.e., the points where the hyper-surface touches the horizontal distribution) 
and characteristic curves.  It turns out that the known (see \cite{Malc}) classification of 
the characteristic points of  the Sub-Riemannian minimal surfaces in the Heisenberg group $H^1$ holds true for the minimal surfaces associated 
to a generic contact $(2,3)$ distribution.  Moreover, we show that in the $(2,3)$ case  the Sub-Riemannian  
minimal surfaces are the integral surfaces of a certain system of ODE in the extended state space. 
In some particular cases the Cauchy problem for this system can be solved explicitly. 
We illustrate our results considering  Sub-Riemannian  minimal surfaces in the Heisenberg 
group and the group of roto-translations.

\end{abstract}
 
\section{Introduction}

In the classical Riemannian geometry minimal surfaces realize the critical points of area functional 
with respect to variations that preserve the boundary of a given domain.
The Sub-Riemannian  minimal surfaces  are the natural  generalization of the classical  ones in Sub-Riemannian manifolds 
known also as the Carnot-Carath\'eodory  spaces.
The notion of a minimal surface in the Sub-Riemannian manifold was introduced  in  \cite{GarN} in the framework of Geometric Measure Theory, 
and then was studied in \cite{GarP},~\cite{FSSC},~\cite{Pau0},~\cite{Pau1},~\cite{Citti},~\cite{Malc},~\cite{Rit}.
The main part of the results of the cited papers are related 
to the Heisenberg group $H^1$, though recently the first steps were done toward the analysis of the  group of roto-translations (\cite{Citti}, \cite{Pau1}).  
A very fruitful geometrical model was recently  proposed in \cite{Malc}. The authors gave a  general geometrical definition of the 
Sub-Riemannian minimal surfaces by means of CR-structures in $3$-dimensional pseudohermitian manifolds and studied in great detail the case of the Heisenberg group $H^1$.
In  our  paper we propose an alternative (with respect to \cite{Malc}) coordinate-free way to 
define Sub-Riemannian minimal surfaces using the tools of Sub-Riemannian geometry\footnote{For the 
detailed exposition on Sub-Riemannian geometry the reader can 
consult \cite {AgrBook},~\cite{AgrExp},~\cite {Bel}.}. 

Our general construction is the following. Denote by  $M$  an $n$-dimensional smooth  manifold  and 
let $\Delta$ be a co-rank 1 smooth vector distribution on it (``horizontal'' distribution). 
Assume that $\Delta$ is endowed with a Riemannian structure, which can be described by fixing an orthonormal 
basis of vector fields $X_1,\dots,X_{n-1}$ on $\Delta$ (see ~\cite{Bel}).
Then we say that $\Delta$ defines a Sub-Riemannian structure on $M$, and $M$ is a Sub-Riemannian manifold.
It turns out that  there is a canonical way to define a volume form $\mu\in\Lambda^n M$ associated to the Sub-Riemannian structure of $M$. 
Moreover, in analogy with the classical Riemannian case one can define the Sub-Riemannian normal of a hyper-surface 
$W\subset M$ as a unite vector field $\nu$ such that 
$$
\intt_\Omega i_{\nu}\mu=\max_{{X\in \Delta}\atop{|X\|_{\Delta}=1}}\intt_\Omega i_X\mu,\qquad \Omega\subset W.
$$
The $n-1$-form $i_{\nu}\mu$ is the Sub-Riemannian analog of the classical area form on $M$.
It is easy to see that this definition correlates perfectly with the classical definition 
of the Riemannian normal and area  since any  Riemannian manifold is  a   Sub-Riemannian manifold
with  $\Delta\equiv TM$.  
 
As in the classical case, one can define the Sub-Riemannian minimal surfaces in $M$ as the critical points
of the  functional associated to the Sub-Riemannian area form.  It turns out that these surfaces 
satisfy  the following  intrinsic equation
$$
(d\circ i_{\nu}\mu)\Big|_W=0.
$$

The described construction opens a wide possibility to study the Sub-Riemannian minimal surfaces associated to generic
Sub-Riemannian structures  of any  dimension. Moreover, it does not require the existence of any additional global structure in $M$.
It worth to mention that in the case of the Heisenberg group, as well as in the case of the 
roto-translational group, our definition coincides with the known ones 
(see \cite{GarP},~\cite{FSSC},~\cite{Pau0}, ~\cite{Citti} and references therein) 

\vspace{10pt}

In the second part of this paper we consider more in detail the case of contact Sub-Riemannian structures, in particular the case 
of  $2$-dimensional distributions in  $3$-dimensional manifolds. In the latter case all information 
related to  the intrinsic geometry of Sub-Riemannian minimal surfaces is encoded in the 
 {\it characteristic curve} $\gamma:\,[0,T]\mapsto W$ such that  $\dot\gamma(t)\in T_{\gamma(t)}W\cap \Delta_{\gamma(t)}$ for all $t\in [0,T]$.  
The tangent of this curve is orthogonal to the Sub-Riemannian normal of $W$ and  it is defined 
(as well as the Sub-Riemannian area and normal $\nu$) away from the {\it characteristic points of $W$}, where $T_q W\equiv \Delta_q$. 
Applying the classical method of characteristics to the minimal surface equation  
we show that the Sub-Riemannian minimal surfaces are actually the integral 
surfaces of a certain system of ODE in the extended space $M\times S^1$. The characteristic 
points are  either the singular points of these surfaces or the singularities of
their projection on the base manifold.  It turns out that  the classification of the characteristic points described in \cite{Malc} (Theorem B) for the Heisenberg case, 
holds true for any  $(2,3)$ contact Sub-Riemannian structure. 

We conclude our analysis comparing the characteristic curves of a  Sub-Riemannian minimal surface 
with the Sub-Riemannian geodesics in $M$. We show that in general the characteristic curves do not coincide with the Sub-Riemannian geodesics, though in some 
particular cases they do, as for example in the case of the Heisenberg group, while in the group of 
roto-translations only a certain class of characteristic curves are Sub-Riemannian geodesics.

The author is grateful to prof. A. Agrachev whose vision of the problem inspired  this work. 
Many thanks also to prof. E. Pagani for stimulating discussions.

\section{Sub-Riemannian minimal surfaces: general construction}

\subsection{Sub-Riemannian structures and associated objects}   

Let $M$ be an $n$-dimensional smooth manifold. Consider a co-rank $1$ vector distribution  $\Delta$ on $M$: 
$$
\Delta=\bigcup_{q\in M}\Delta_q,\quad \Delta_q\subset T_qM,\quad q\in M.
$$
By definition, the Sub-Riemannian structure on $M$  is a pair $(\Delta, \langle\cdot,\cdot\rangle_{\Delta})$, where
$\langle\cdot,\cdot\rangle_{\Delta}$ denotes a  smooth family of Euclidean inner products on $\Delta$.
In what follows we will call $\Delta$ {\it the horizontal distribution} and  keep 
the same notation $\Delta$ both for the vector distribution and for the associated Sub-Riemannian structure.

Let $X_i$, $i=1,\dots, n-1$ be a horizontal orthonormal basis on $\Delta$:
$$
\Delta_q={\rm span}\{X_1(q),\dots, X_{n-1}(q)\},\qquad q\in M,
$$
$$
\la X_i(q),X_j(q)\ra_{\Delta}=\delta_{ij},\qquad q\in M,\;i,j=1,\dots, n-1.
$$
By $\Theta\in \Lambda^{n-1}\Delta$ we will denote the corresponding Euclidean volume form on $\Delta$.

In what follows we will assume that $\Delta$ is {\it bracket-generating} on $M$. In the present case this means that
$$
{\rm span}\{X_i(q), [X_i, X_j](q),\;i,j=1,\dots,n-1,\;q\in M\}=T_q M.
$$
Hereafter the square brackets denote the Lie brackets of vector fields.
If $\Delta$ is bracket-generating, then by the Frobenius theorem it is completely non-holonomic, i.e.,
there are no invariant sub-manifolds of  $M$ such that their tangent spaces coincides with 
$\Delta$ at any point.

\vspace{10 pt}

There is an alternative way to define the  distribution $\Delta$  
as the kernel of some differential $1$-form. Let  $\omega\in \Lambda^1M$ be such a form : 
$$
\Delta_q={\rm Ker}\, \omega_q=\{v\in T_q M:\;\omega_q(v)=0\},\quad q\in M.
$$
In general, the form $\omega$  is defined up to a multiplication by a non-zero scalar function. 
It  is easy to check that $\Delta$ is bracket-generating at $q\in M$ 
if and only if $d_q\omega\ne 0$.

By the standard construction the Riemannian structure on $\Delta$ can be extended to 
the spaces of forms $\Lambda^k \Delta$, $k\le n-1$.
In particular, for any  $2$-form $\sigma$ we set
$$
\|\sigma_q\|_{\Delta}=\left(\sum_{i,j=1\atop i<j}^{n-1} \sigma_q(X_i(q), X_j(q))^2\right)^{\frac{1}{2}},
$$
$\{X_i(q)\}_{i=1}^{n-1}$, as before, being  an orthonormal horizontal basis of $\Delta_q$.
Now we  can fix  the choice of the form $\omega$ defined above  by setting  
\BE\label{omega}
\omega_q(\Delta_q)=0,\qquad \|d_q\omega\|_{\Delta}=1,\qquad \forall q\in M.
\EE
We will call the $1$-form satisfying (\ref{omega}) {\it the  canonical $1$-form associated to $\Delta$}.
In the fixed horizontal orthonormal basis $\{X_i(q)\}_{i=1}^{n-1}\in \Delta_q$ equations (\ref{omega}) become
\BE\label{omega1}
\omega_q(X_i(q))=0,\qquad \sum_{i,j=1\atop i<j}^{n-1} d_q\omega(X_i(q), X_j(q))^2=1,\qquad i=1,\dots, n-1.
\EE
In worth to note that  the canonical $1$-form, $\omega$ satisfying (\ref{omega1})
is defined up to a sign and does not depend on the choice of the horizontal basis.
In local coordinates in $M$ the components of $\omega$ can be expressed in terms of  the coordinates of the vector 
fields $X_i$ and their first derivatives, since due to the Cartan formula 
$$
d\omega(X,Y)=X \omega(Y)-Y \omega(X)-\omega([X,Y]),
$$ 
and hence
$$
\|d_q\omega\|_\Delta^2=\sum_{i,j=1\atop i<j}^{n-1} d_q\omega(X_i(q), X_j(q))^2=\sum_{i,j=1\atop i<j}^{n-1}\omega_q([X_i,X_j](q))^2.
$$

Once the orientation in $M$ if fixed by choosing the sing of $\omega$, the following volume form
$$
\mu=\Theta\wedge \omega
$$
is uniquely defined. We will call this volume form {\it the canonical volume form associated to $\Delta$}. 
The canonical volume form  $\mu$ is a ``global'' object  in $M$, though it is intrinsically defined by the Sub-Riemannian structure on $\Delta$.

\subsection{Horizontal area form}

Let $W\subset M$, ${\rm dim} W=n-1$ be a smooth hyper-surface in $M$ and let $\Omega\subset W$ be an open domain.
For simplicity we assume that the  vector field $X\in TM$ is transversal to $W$, though this assumption is not restrictive. Consider the flow  generated by $X$ in $M$ :
$$
\Pi^X:\qquad [0,\varepsilon]\times \Omega\mapsto M,
$$
$$
\Pi^X(t,q)=e^{t\,X}(q),\quad q\in M.
$$
Denote by 
\BE
\Pi_{(\eps, \Omega)}^X=\left\{e^{tX}(q),\;q\in\Omega,\;t\in[0,\varepsilon]\right\} 
\EE
the cylinder formed by the images of $\Omega$ translated along the integral curves of $X$ 
parametrized by $t\in[0,\eps]$. Clearly, $\Pi_{(0, \Omega)}^X\equiv\Omega$.
By definition,
$$
Vol(\Pi_{(\eps, \Omega)}^X)=\int\limits_{\Pi^X_{(\eps, \Omega)}}\mu =\int\limits_{[0,\varepsilon]\times\Omega}({\Pi^X})^*\mu,
$$  
where $({\Pi^X})^*$ denotes the pull-back map associated to $\Pi^X$ and $\mu$ is the canonical volume form defined above\footnote{Here  
we use the canonical volume form associated to $\Delta$, though the whole construction works for any volume form in $M$.}.

\begin{defi} The following quantity
\BE\label{areaform}
A_\Delta(\Omega)=\max_{{X\in \Delta}\atop{\|X\|_{\Delta}=1}}\;\lim_{\varepsilon\to 0}\frac{Vol(\Pi^X_{(\eps,\Omega)})}{\eps}
\EE
is called the Sub-Riemannian  (or {\it horizontal}) area  of the domain $\Omega$ associated to $\Delta$.\end{defi}

\noindent{\bf Remark.} The horizontal area defined by (\ref{areaform}) is nothing but the generalization of the classical notion of the Euclidean area:
 it defines the area of the base of a cylinder as the ratio of its volume and height.

Let us find a more convenient expression for (\ref{areaform}). First of all we observe that since 
$({\Pi^X})^*\mu$ is a form of maximal rank $n$ on $M$ we have  $dt\wedge({\Pi^X})^*\mu=0$. Hence
$$
0=i_{\dd_t}\left(dt\wedge({\Pi^X})^*\mu\right)=i_{\dd_t}dt\wedge({\Pi^X})^*\mu-dt\wedge i_{\dd_t}({\Pi^X})^*\mu,
$$
i.e.,
$$
({\Pi^X})^*\mu=dt\wedge i_{\dd_t}({\Pi^X})^*\mu.
$$
Taking into account that $\Pi^X_*\dd_t=X$ we obtain
$$
A_\Delta(\Omega)=\max_{{X\in \Delta}\atop{\|X\|_{\Delta}=1}}\;\lim_{\varepsilon\to 0}\frac{Vol(\Pi^X_{(\eps,\Omega)})}{\eps}=
$$
$$
= \max_{{X\in \Delta}\atop{\|X\|_{\Delta}=1}}\;\frac{\dd}{\dd \eps}\Big|_{\eps=0}\intt_{\Pi^X_{(\eps,\Omega)}}\mu
=\max_{{X\in \Delta}\atop{\|X\|_{\Delta}=1}}\;\frac{\dd}{\dd \eps}\Big|_{\eps=0} \intt_{[0,\eps]\times\Omega}dt\wedge i_{\dd_t}({\Pi^X})^*\mu=
$$
$$
= \max_{{X\in \Delta}\atop{\|X\|_{\Delta}=1}}\;\frac{\dd}{\dd \eps}\Big|_{\eps=0}\intt_0^{\eps}\Big(\intt_{\Pi^X_{(t,\Omega)}}i_X\mu\Big) dt=
 \max_{{X\in \Delta}\atop{\|X\|_{\Delta}=1}}\intt_\Omega i_X\mu.
$$

\begin{defi}
The horizontal unite  vector field $\nu\in\Delta$, $\|\nu\|_{\Delta}=1$ such that
$$
\intt_\Omega i_{\nu}\mu=\max_{{X\in \Delta}\atop{|X\|_{\Delta}=1}}\intt_\Omega i_X\mu
$$
is called the Sub-Riemannian  or horizontal normal of $\Omega\subset W$. 
The $(n-1)$-form $i_{\nu}\mu$ is called the Sub-Riemannian  or horizontal area form associated to $\Delta$.
\end{defi}

\vspace{10pt}

In general,  the given  definition of the Sub-Riemannian normal $\nu$ does not require the existence of any 
global structure in $M$ (for instance, one does not need a Riemannian structure in $M$).
Nevertheless, if $M$ is a Riemannian manifold whose Riemannian structure is compatible with the Sub-Riemannian
structure on $\Delta$, i.e., if the inner product $\la\cdot, \cdot\ra$ on $TM$ satisfies $\la\cdot,\cdot\ra_{\Delta}=\la\cdot, \cdot\ra\big|_\Delta$, 
then it is easy to see  that the Sub-Riemannian normal $\nu$ is  nothing but the projection on $\Delta$ of the Riemannian 
unit normal $N$ of $W$, normalized w.r.t. $\|\cdot\|_{\Delta}$. This is the consequence of the following  relation
$$
\int\limits_{\Omega} i_{X} \mu =\int\limits_{\Omega}\langle X,N\rangle i_N\,\mu,\qquad \forall  X\in Vec(M).
$$
Thus if $X_1,\dots,X_{n-1}\in \Delta$ is an orthonormal horizontal basis of $\Delta$,
then 
\BE\label{srnorm}
\nu=\sum\limits_{i=1}^{n-1}\nu_i X_i\,,\qquad \nu_i=\frac{\langle N,X_i\rangle}
{\sqrt{\langle N,X_1\rangle ^2+\dots+\langle N,X_{n-1}\rangle ^2}},
\EE
and the  horizontal area form reads
$$
i_{\nu}\mu=\langle \nu, N \rangle i_N \,\mu=
\sqrt{\langle N,X_1\rangle ^2+\dots+\langle N,X_{n-1}\rangle ^2}\,i_N\,\mu.
$$
Now let the hyper-surface $W$ be defined as a level set  of a smooth, let us say $C^2$, function:
$$
W=\left\{q\in M:\quad F(q)=const,\;F\in C^2(M), \;d_q F\ne 0\right\}.
$$
Let $X \equiv X_n$ be a vector field transversal to $W$ and such that $\{X_1(q),\dots,X_n(q)\}$ 
form an orthonormal basis of $T_q M$ at $q\in M$. 
Then 
$$
N(q)=D_0^{-1}\sum\limits_{i=1}^n X_i F(q)\,X_i(q),\qquad D_0=\left(\sum\limits_{i=1}^n X_i F(q)^2\right)^{1/2}
$$
and
\BE\label{nu}
\nu(q)=D_1^{-1}\sum\limits_{i=1}^{n-1} X_i F(q)\, X_i(q) ,\qquad D_1=\left(\sum\limits_{i=1}^{n-1} X_i F(q)^2\right)^{1/2}.
\EE
Here  $X_iF$ denotes the directional derivative of $F$ along the vector field $X_i$.

\subsection{Sub-Riemannian minimal surfaces}

Now let us  compute the variation of the horizontal area $A_\Delta(\cdot)$.
Assume that  $\Omega\subset W$ is a bounded domain  and let $V\in Vec(M)$ be such that $V\big|_{{\partial}\Omega}=0$.
Consider a one-parametric family of hyper-surfaces generated by the vector field $V$
$$
\Omega^t= e^{t V}\Omega,\qquad \Omega^0\equiv\Omega,
$$
and denote by $\nu^t$ the horizontal unit normals to $\Omega^t$.
We have
$$
A_{\Delta}(\Omega^t)=\int\limits_{e^{t V}\Omega} i_{\nu^t}\ \mu=\int\limits_{\Omega}(e^{t V})^* i_{\nu^t}\,\mu
=\int\limits_{\Omega}e^{t L_V} i_{\nu^t}\,\mu.
$$
Further, 
\BE\label{int1}
\frac{\partial}{\partial t} \Big|_{t=0} A_{\Delta}(\Omega^t)=\int\limits_{\Omega}L_V i_{\nu}\,\mu+
\int\limits_{\Omega}i_{\frac{\partial \nu^t}{\partial t}\big|_{t=0}}\mu.
\EE
It is not hard to show that the second integral in (\ref{int1}) vanishes, 
because the horizontal vector field  $\frac{\dd \nu^t}{\dd t}\big|_{t=0}$ is tangent to $\Omega$.
Indeed, at a generic (non-characteristic) point $q\in \Omega$ we have $\nu(q)\notin T_q \Omega$ and 
${\rm dim}\Delta_q\cap T_q \Omega=n-2$. On the other hand,  differentiating 
the equality $\la \nu^t,\nu^t\ra_{\Delta}=1$ we get 
\BE\label{unity}
\langle \frac{\partial \nu^t}{\partial t}\Big|_{t=0},\nu\rangle_{\Delta}=0,
\EE
and hence  $\frac{\dd \nu^t}{\dd t}\big|_{t=0}(q)\in T_q \Omega$.

Further, using Cartan's formula we transform the first part of (\ref{int1}):
$$
\int\limits_{\Omega}L_V i_{\nu}\,\mu=\int\limits_{\Omega}(i_V\circ d+d\circ i_V) i_{\nu}\,\mu=
\int\limits_{\Omega}(i_V\circ d \circ i_{\nu})\mu+\int\limits_{\Omega}(d\circ i_V\circ i_{\nu})\mu.
$$
Applying the Stokes theorem to the second integral we see that it vanishes:
$$
\int\limits_{\Omega}(d\circ i_V\circ i_{\nu})\mu=\int\limits_{\partial \Omega}(i_V\circ i_{\nu})\mu=0
$$
provided $V\big|_{\partial \Omega}=0$ and $\dd \Omega$ is sufficiently regular.
Thus, 
$$
\frac{\partial}{\partial t} \Big|_{t=0} A_{\Delta}(\Omega^t)=\int\limits_{\Omega}i_V\circ(d \circ i_{\nu}\mu).
$$

\begin{defi} We  say that the hyper-surface $W$ is a minimal surface w.r.t. the Sub-Riemannian structure $\Delta$ (or just $\Delta$-minimal) iff     
\BE\label{mineq}
( d\circ i_{\nu}\mu)\Big|_{W}=0.
\EE
\end{defi}

We remark that the minimality of a hyper-surface does not depend on the chosen orientation in $M$.
It is also easy to see that the whole construction can be easily generalized for the case of vector distributions of co-rank greater than $1$.

\subsection{Canonical form of the minimal surface equation in contact Sub-Riemannian  manifolds}

Let  $n=2m+1$ and assume that $\Delta$ is a contact distribution, 
i.e., the $2m+1$-form $(d\omega)^m\wedge \omega$ is non-degenerate. Then we say that $M$  is a  contact Sub-Riemannian manifold.
In the contact case there exists a unique vector filed $X\in TM$ such that
\BE\label{Reeb}
\omega_q(X(q))=1, \qquad d_q\omega(V, X(q))=0,\qquad \forall V\in \Delta_q.
\EE
Such a vector field is called {\it the Reeb vector field} associates to the contact form $\omega$.
Using this vector field we can extend the Sub-Riemannian structure on $\Delta$ to the whole $TM$. The resulting Riemannian structure in $M$ 
is by definition compatible with $\Delta$. The basis of vector fields $\{X_1,\dots, X_{2m},X\}$ is  then  a canonical basis 
associated to the contact Sub-Riemannian structure $\Delta$.

\vspace{10pt}

Set  $X_{2m+1}\equiv X$. Denote by $c_{ij}^k\in C^\infty(M)$ {\it the  structural constants} of the frame $\{X_i\}_{i=1}^{2m+1}$ 
\BE\label{se}
[X_i, X_j]=-\sum\limits_{k=1}^{2m+1} c_{ij}^k X_k.
\EE
Let $\{\theta_i\}_{i=1}^{2m+1}$be  the basis of $1$-forms dual to $\{X_i\}_{i=1}^{2m+1}$. 
Clearly, $\theta_{2m+1}\equiv\omega$ and the canonical volume form is
$$
\mu=\theta_1\wedge\dots\wedge\theta_{2m+1}.
$$
From the Cartan formula it follows that
\BE\label{Cart}
d\theta_k=\sum_{i,j=1\atop i<j}^{2m+1}c_{ij}^k\theta_i\wedge\theta_j,\qquad k=1,\dots, 2m+1. 
\EE
 
Let us now derive the canonical form of the minimal surface equation (\ref{mineq}) in contact Sub-Riemannian manifolds.
First we calculate the interior product of $\nu=\sum\limits_{i=1}^{2m} \nu_i X_i$ with the 
canonical volume form:
$$
i_\nu\,\mu=\left(\sum_{k=1}^{2m}(-1)^{k+1}\nu_k \,\theta_1\wedge\dots\wedge\widehat{\theta_k}\wedge\dots\wedge \theta_{2m}\right)\wedge\theta_{m+1}=\Xi\wedge\theta_{2m+1}.
$$
Here $\widehat{\theta_k}$ denotes the omitted element in the wedge product and 
$$\Xi=\sum_{k=1}^{2m}(-1)^{k+1}\nu_k \,\theta_1\wedge\dots\wedge\widehat{\theta_k}\wedge\dots\wedge \theta_{2m}.$$
Further, 
$$
d\,i_\nu\,\mu=d \Xi\wedge\theta_{2m+1}-\Xi\wedge d\theta_{2m+1}.
$$
Recalling now that $d\nu_k=\sum\limits_{i=1}^{2m+1}\,X_i(\nu_k)\theta_i$, we obtain
$$
d\Xi\wedge\theta_{2m+1}=\sum_{k=1}^{2m}(-1)^{k+1}\left(d\nu_k\wedge \theta_1\wedge\dots\wedge\widehat{\theta_k}\wedge\dots\wedge \theta_{2m}+\right.
$$
$$
\left.+\nu_k\, d(\theta_1\wedge\dots\wedge\widehat{\theta_k}\wedge\dots\wedge \theta_{2m})\right)\wedge\theta_{2m+1}=
\left(\sum_{k=1}^{2m}X_k(\nu_k)+\sum_{j=1}^{2m}\nu_k c_{kj}^j\right)\mu.
$$
On the other hand, 
$$
\Xi\wedge d\theta_{2m+1}=\Xi\wedge \sum_{i,j=1\atop i<j}^{2m+1}c_{ij}^{2m+1}\theta_i\wedge\theta_j=-\left(\sum_{k=1}^{2m}\nu_k c_{k 2m+1}^{2m+1}\right)\mu.
$$
Summing up we obtain the following equation:
\BE\label{mineqN}
\left.\left[{\rm div}^{\Delta}\nu+\sum\limits_{i=1}^{2m} \nu_i(q)\left(\sum\limits_{j=1}^{2m+1} c_{ij}^j\right)\right]\right|_W=0. 
\EE

The left-hand side of  (\ref{mineqN}) is called {\it the Sub-Riemannian mean curvature} of the hyper-surface $W$, while its first term
$$
{\rm div}^{\Delta}\nu=\sum_{i=1}^{2m} X_i(\nu_i)
$$
is called {\it the horizontal divergence} of the Sub-Riemannian normal $\nu$. 
Equation (\ref{mineqN}) is the canonical equation of Sub-Riemannian minimal surfaces in a contact Sub-Riemannian manifold.
In the rest of the present paper we will try to analyze it in the less-dimensional case of $m=1$.

\vspace{10pt}

\noindent{\bf Remark} In general the Sub-Riemannian structures are not equivalent to the CR-structures, which were 
used in \cite{Malc} and the successive publications by other authors, and consequently in general equation (\ref{mineqN})
is different from its analog obtained  in \cite{Malc} for $2$-dimensional minimal 
surfaces in  $3$-dimensional contact CR manifolds. Nevertheless, in some particular cases,  like the 
Heisenberg group and the group of roto-translations, both models produce the same  result.

\label{examp1}\begin{example}(The Heisenberg distribution) {\rm
Let $M=\real^{2m+1}$ and denote by  $(x_1,\dots, x_{2m},t)=q$ the Cartesian coordinates in $M$.
Let  $\Delta$ be such that $\Delta_q={\rm span}\{X_i(q)\}_{i=1}^{2m}$, $q\in M$, where
\BE\label{heis}
X_i(q)=\dd_{x_i}+\frac{x_{i+m}}{2}\dd_t,
\EE
$$
X_{i+m}(q)=\dd_{ x_{i+m}}-\frac{x_i}{2}\dd_t,\quad i=1,\dots,m.
$$
The vector distribution $\Delta$ is characterized by the following commutative relations:
\BE\label{deltaH}
[X_i,X_j]=0,\quad {\rm for}\quad j\ne i+m,\qquad [X_i,X_{i+m}]=-\dd _t,
\EE
and therefore it is a co-rank $1$ bracket-generating distribution.
The vector fields $X_i$, $i=1,\dots,2n$, generates the so-called {\it Heisenberg Lie algebra} on $\real^{2m+1}$.
In what follows we will call the vector  distributions which satisfy the commutative relations (\ref{deltaH}) 
{\it the Heisenberg distribution} and denote it by $\Delta^{H^m}$. The space
$\real^{2m+1}$ endowed with the structure of this distribution is called the Heisenberg group $H^{m}$.

By solving (\ref{omega1}) we find the canonical $1$-form $\omega$: 
\BE\label{omegaH}
\omega=\pm\frac{1}{\sqrt{m}}(dt -\frac{1}{2}\sum_{i=1}^m(x_{i+m}\,dx_i-x_i\,dx_{i+m})), 
\EE
and correspondingly the Reeb vector field $X=\pm \sqrt m \dd_t$.  Clearly $\omega$ is a contact form since 
$(d\omega)^m\wedge d\omega=\pm\frac{1}{m^m}\bigwedge\limits_{i=1}^{2m}dx_i\wedge dt$ in non-degenerate. 
The only non-zero structural constants of the canonical frame are   $c_{i i+m}^{2m+1}=\pm \frac{1}{\sqrt m}$, $i=1,\dots,m$.
Due to the hight degeneracy of the Sub-Riemannian structure the canonical 
minimal surface equation takes a very simple form:   
$$
{\rm div}^{\Delta^{H^m}}\nu\;\Big|_W=0.
$$
This is the well known minimal surface equation in the Heisenberg group (see \cite{GarP}, \cite{FSSC}, \cite{Malc}, \cite{Rit}, etc.)
}\end{example}

\section{Sub-Riemannian minimal surfaces for $(2,3)$ contact vector distributions}

In this section we analyze the case of a contact distribution $\Delta$ of rank $2$ in the 
$3$-dimensional manifold $M$. In this case the intrinsic information about the geometry of the 
$\Delta$-minimal surface  $W$ is encoded in the so-called {\it characteristic curves} of $W$, 
which can be defined as the leaves of the one-dimensional foliation 
$TW\cap \Delta$. The singular points of the  characteristic curves are called {\it the characteristic points}. 
At these points $\Delta$ is tangent to $W$ and hence the Sub-Riemannian normal 
(as well the horizontal area form) is not defined.  

In the case of $(2,3)$ contact distributions the characteristic curve, being a one-dimensional sub-manifold, 
has no intrinsic invariants. However, one can extract some information about 
the global geometry of the $\Delta$-minimal
surfaces by analyzing the type of its characteristic points.

\vspace{10pt}
 
Let $n=3$ and assume that $\Delta$ is such that $\Delta_q={\rm span}\{X_1(q),X_2(q)\}$, $q\in M$.
Set $X_3\equiv X$, where $X$ is the Reeb vector field associated to $\Delta$, and denote by  $c_{ij}^k$ the structural constant of the canonical frame 
$\{X_i\}_{i=1}^3$.
By definition, $c_{ij}^k=-c_{ji}^k$.
Moreover, (\ref{Reeb}) and (\ref{se}) imply
\BE\label{se3}
c_{12}^3=1,\qquad c_{13}^3=c_{23}^3=0.
\EE
More symmetry  relations of the structural constants can be obtained from the 
Jacobi identity
$$
[X_1,[X_2,X_3]]+[X_3,[X_1,X_2]]+[X_2,[X_3,X_1]]=0.
$$
In particular, if $M$ is a Lie group, the structural constants do not depend on the points of the base manifold $M$, and
the Jacobi identity is equivalent to the following relations:
\BE\label{jac}
c_{13}^1+c_{23}^2=0,\qquad
c_{12}^1c_{13}^1+c_{12}^2c_{23}^1=0,\qquad
c_{12}^1c_{13}^2+c_{12}^2c_{23}^2=0.
\EE

\vspace{10pt}

Let $\nu\in \Delta$ be a horizontal normal of a regular hyper-surface $W\subset M$. 
Taking into account (\ref{se3}), we  write the $\Delta$-minimal surface equation
at non-characteristic points:
\BE\label{mineq1}
\Big({\rm div}^{\Delta}\nu+\nu_1 c_{12}^2 - \nu_2 c_{12}^1\Big)\Big|_W=0.
\EE
If $W$ is a level set of some smooth function $F$, using  (\ref{nu}), we obtain
\BE\label{mineqF}
\Big[\Big(X_1^2F\,(X_2 F)^2+X_2^2 F\,(X_1F)^2-X_1F\, X_2 F\,(X_1\circ X_2+X_2\circ X_1)F\Big)D_1^{-3}+
\EE
$$
+\Big(c_{12}^2 X_1F-c_{12}^1X_2F\Big)D_1^{-1}\Big]\Big|_W=0,\qquad D_1=\sqrt{X_1F^2+X_2F^2}. 
$$
The last equation is a highly degenerate PDE.  
Many non-trivial solutions of this equation are known for $\Delta^{H^1}$, the interested reader can consult
in \cite{Malc} and other papers, cited in the Bibliography. Another important for applications case is the distribution, 
which corresponds to another Lie Group, the so-called {\it roto-translational group} $e^2$. 

\label{examp2}\begin{example}{\rm The Lie group $e^2$ can be realized as $\real^2\times {\mathbb S}^1$ 
with  local coordinates $(x,y,z)$. The Lie algebra corresponding to this group is generated by  
vector fields 
$$
X_1=\cos z \dd_x+\sin z\dd_y,\qquad  X_2=\dd_z.
$$
It is easy to check that the  horizontal distribution  $\Delta^{e^2}$ with sections $\Delta^{e^2}_q={\rm span}\{X_1(q),X_2(q)\}$, $q\in M$, is contact, 
the corresponding canonical $1$-form is $\omega=\pm(\sin z dx -\cos z dy)$. 
The Reeb vector field coincides with the Lie bracket $[X_1,X_2]$ (up to the sign) 
and the only non-zero structural constants are 
$c_{12}^3=c_{23}^1=\pm 1$. The following  surfaces are $\Delta^{e^2}$-minimal surfaces 
(away from the characteristic points):\\
a). $y=x+ B(\sin z+\cos z)+C,\qquad B,C={\rm const}$;\\
b). $A x+ B \sin z= C,\qquad A,B,C={\rm const}$;\\
c). $x \cos z+ y \sin z=0$.

}\end{example}

\subsection{Structure of characteristic points}

In \cite{Malc} the authors showed  (see Theorem B) that in the case  
of the Heisenberg distribution $\Delta^{H^1}$ the characteristic points of the corresponding minimal surfaces are 
either isolated of index $+1$ or contained in a $C^1$ curve. Such curves are called singular. The characteristic curves  keep go straight after they cross a singular curve.
All these facts are of local nature, and it turns out that they hold true for Sub-Riemannian minimal surfaces in generic contact Sub-Riemannian manifolds of dimension $3$.
Basically all local arguments used in \cite{Malc} can be directly applied in the general case  modulo a  suitable choice of local coordinates 
in $M$ in the small neighborhood of a characteristic point.

\vspace{10pt}

Let $F\in C^2(M)$ and let $W$ be a $\Delta$-minimal smooth surface in $M$ defined as a level set of $F$, i.e., assume $F$ satisfies (\ref{mineqF}) away from the points
$\hat q\in W$ where $X_1F(\hat q)=X_2F(\hat q)=0$. 

Let $\hat q$ be a characteristic point of $W$. Since  $d_{\hat q}F\ne 0$ one can choose the local coordinates in $M$ in a small neighborhood 
${\mathcal O}_{\hat q}=\{q=(x,y,z)\in \real^3\}$ in such a way that $W=\{q;\;F(q)\equiv z=0, \;q\in{\mathcal O}_{\hat q}\}$. 
Since $\hat q$ is a characteristic point $\Delta_{\hat q}=T_{\hat q}W=\{(x,y,0)\}$. 
Then
$$
X_1(x,y,z)=\dd_x+a_1(x,y,z)\dd_z,$$
\BE\label{coord}
X_2(x,y,z)=x a_3(x,y,z)\dd_x +\dd_y+a_2(x,y,z)\dd_z\,,
\EE
$$
a_i(0,0,0)=0\qquad{\rm for}\quad i=1,2,
$$
where $a_1$, $a_2$ and $a_3$  are some functions, which define the Sub-Riemannian structure on $\Delta$.
In particular, since $\Delta$ is bracket generating at any point
$$
{\rm rank}\{X_1(\hat q),X_2(\hat q),[X_1,X_2](\hat q)\}=3,
$$ 
and so
\BE\label{det}
{\rm det}\left(\begin{array}{ccc}
1&0&0\\
0&1&0\\
a_3(0)&0&\frac{\dd a_2}{\dd x}(0)-\frac{\dd a_1}{\dd y}(0)\end{array}\right)=\frac{\dd a_2}{\dd x}(0)-\frac{\dd a_1}{\dd y}(0)\ne 0.
\EE
Denote
$$
A=\left(\begin{array}{cc}
\frac{\dd a_1}{\dd x}&\frac{\dd a_1}{\dd y}\\
\frac{\dd a_2}{\dd x}&\frac{\dd a_2}{\dd y}\end{array}\right)(0)=
\left(\begin{array}{cc}
a&b\\
c&d\end{array}\right).
$$
The matrix $A$ plays a key role on the further analysis.  First of all, condition (\ref{det}) means that  $b$ and $c$ cannot be zero simultaneously. 
Thus  $A$ is not a $0$-matrix and  ${\rm trace}A\ne 0$. Moreover, the implicit function theorem implies that  ${\rm det A}=0$ 
if and only if the point $\hat q$ is not isolated  and contained in a $C^1$ curve called singular. In particular,  if $\hat q$ is isolated, then 
${\rm det A}\ne 0$.  If $\hat q$ is not isolated, the characteristic curves keep go straight through a singular curve.
The proof of these facts is rather technical and it repeats  the proof of the  
analog result in the Heisenberg group ( see \cite{Malc}) with a few changes due to 
the use of the new curvilinear coordinates $(x,y,z)$. We omit this proof here, an interested reader can  find it in the original paper.
Let us just show explicitly that the index of an isolated characteristic point of a $\Delta$-minimal surface is equal to $+1$ and it not affected by the difference of 
the given Sub-Riemannian structure on $\Delta$  from the Heisenberg one.

By definition, the characteristic point $\hat q$ is the singular point of the vector field $\nu_0=X_1F\,X_1+X_2F\,X_2$ and
${\rm ind}(\hat q)={\rm sgn\;det} A$. Denote $H= X_1^2 F (X_2 F)^2+X_2^2 F (X_1F)^2-X_1F\,X_2 F\,(X_1\circ X_2+X_2\circ X_1)F$. Then equation (\ref{mineqF}) becomes
\BE\label{mineqF1}
\frac{H+(X_1 F^2+X_2 F^2)(c_{12}^2 X_1 F-c_{12}^1 X_2 F)}{D_1^3}=0,
\EE
Let $q=\hat q+\delta q$, where $\delta q=(\delta x,\delta y,0)$ and  $r=\sqrt{\delta x^2+\delta y^2}$.
By (\ref{coord}), we have
$$
X_1 F(q)=a\delta x+b\delta y +o(r),\qquad X_2 F(\hat q)=c\delta x+d\delta y+o(r),
$$
$$
X_1^2 F(q)=a,\qquad X_2^2 F(q)=d,\qquad(X_1\circ X_2+X_2\circ X_1)F(q)=c+b,
$$
$$
D_1^2(q)=(X_1 F^2+X_2 F^2)(q)=
$$
$$
=(a^2+c^2)\delta x^2+(b^2+d^2)\delta y^2+2(ab+cd)\delta x\delta y+o(r^2).
$$
Substituting these expressions into (\ref{mineqF1}) we obtain
$$
\big(H+(X_1 F^2+X_2 F^2)(c_{12}^2 X_1 F-c_{12}^1 X_2 F)\big)(q)=
$$
$$
=a\big(c\delta x+d\delta y+o(r)\big)^2+d\big(a\delta x+b\delta y+o(r)\big)^2-
$$
$$
-(c+b)\big(a\delta x+b\delta y+o(r))(c\delta x+d\delta y+o(r)\big)+
$$
$$
+\big(c_{12}^2(\hat q +\delta q)(a\delta x+b\delta y+o(r))-
c_{12}^1(\hat q+\delta q)(c\delta x+d\delta y+o(r))\big)\times
$$
$$
\times\big((a^2+c^2)\delta x^2+(b^2+d^2)\delta y^2+2(ab+cd)\delta x\delta y+o(r^3)\big)=
$$
$$
=(ad-bc)(a\delta x^2+d\delta y^2+(c+b)\delta x\delta y)+o(r^2).
$$
We see that the structural constants $c_{12}^1$ and $c_{12}^2$ enters into the play together 
with the higher order terms and do not affect the type of characteristic points.
Further, since $\hat q$ is isolated ${\rm det} A=ad-bc\ne 0$. Moreover, by (\ref{mineqF1}), we have 
$$
\frac{a\delta x^2+d\delta y^2+(c+b)\delta x\delta y+o(r^2)}{((a^2+c^2)\delta x^2+(b^2+d^2)\delta y^2+2(ab+cd)\delta x\delta y+o(r^2))^{3/2}}=0.
$$
Observe that if we take $\delta x=0$, in order to satisfy the last equation we should necessarily have $d=0$, analogously $\delta y=0$ forces $a=0$. 
This implies $b+c=0$ and hence $b=-c$. Therefore ${\rm det A}=-bc=c^2$ and  hence ${\rm ind}(\hat q)=+1$.

So, locally the  $\Delta^{H^1}$-minimal surfaces give a good approximation of the structure of generic  $\Delta$- minimal surfaces. But globally
this is not true.  We well show this difference in the next subsection.

\subsection{Characteristic curves}

Recall that the curves formed by the intersection of the distribution $\Delta$ with $W$ are called {\it the characteristic curves} of the $\Delta$-minimal
surface $W$.  These curves  are  the integral curves of {\it the  characteristic  vector field} $e\in \Delta$ 
such that $\la e,\nu\ra_{\Delta}=0$ and $\|e\|_\Delta=1$.
Thus $e=e_1 X_1+e_2 X_2$ and we fix the orientation on $\Delta$ by setting  $e_1=\nu_2$, $e_2=-\nu_1$.
Since $\|e\|_{\Delta}=1$ one can  introduce an auxiliary  parameter $\phi\in S^1$ 
such that $\cos \phi=e_1$, $\sin\phi=e_2$, and
$$
e\equiv e^{\phi}=cos \phi X_1+\sin\phi X_2.
$$
Here we use the  upper-index $\phi$ to stress out the dependence of the vector field $e$ on $\phi$.
Now equation (\ref{mineq1}) can be rewritten as follows:
\BE\label{mineq3}
-\sin\phi\, X_2 \phi-\cos\phi\, X_1\phi=\cos\phi\, c_{12}^1+\sin\phi\, c_{12}^2.
\EE
Equation (\ref{mineq3}) is a quasilinear PDE and one can apply the classical method of characteristics to find its solutions.
Indeed, let $s\mapsto (q_1(s), q_2(s), q_3(s))$ be a smooth (at least $C^1$) 
curve in $M$. Along this curve $\dot \phi=\sum\limits_{i=1}^3\frac{\dd \phi}{\dd q_i}\dot q_i$ with $\dot{\;\;}=\frac{d}{ds}$. Then (\ref{mineq1}) is equivalent to 
the following system of the first order ODE:
\BE\label{sys}
\left\{\begin{array}{ccl}
\dot q&=&e^{\phi}(q)\\
\dot \phi&=&-\cos\phi\, c_{12}^1(q)-\sin\phi\, c_{12}^2(q)
\end{array}\right.
\EE
This shows that the $\Delta$-minimal surface $W$ in nothing but the projection on $M$ of the integral surface of the system (\ref{sys}) in the extended  
manifold $M_0=M\times S^1=\{(q,\phi):\,q\in M,\,\phi\in S^1\}$. Thus, at least locally, one can find  solutions of the $\Delta$-minimal surface 
equation by solving the Cauchy problem for the system (\ref{sys}). The characteristic points  of $W$, being the points where (\ref{mineq}) 
is not defined, are either the singular points of the surface (\ref{sys}) in the extended space, 
or the singular points of the projection $\pi:M_0\mapsto M$, for instance, the points where the different characteristics meet each other.

\vspace{10pt}

In a particular but important for the applications case of Lie groups  system (\ref{sys}) can be integrated explicitly, at least formally.
Indeed, if $c_{12}^1$ and $c_{12}^2$ are constant, the second equation of (\ref{sys}) can be integrated separately and then the obtained function $\phi(t)$ 
can be used to integrate the first equation of (\ref{sys}).     
Moreover, if
\BE\label{ruledcond}
c_{12}^1(q)=c_{12}^2(q)=0\qquad \forall q\in M,
\EE
then $\phi$ is constant along any characteristic curve. The corresponding minimal surface is a kind of ruled surface, whose rulings are the characteristic curves
that are not straight lines in general.

\Grarray{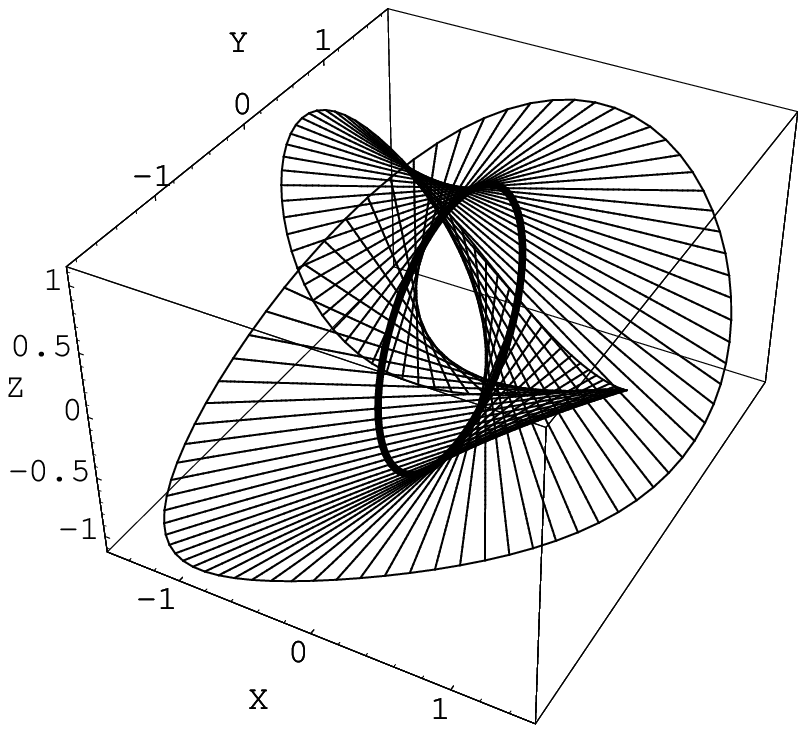}{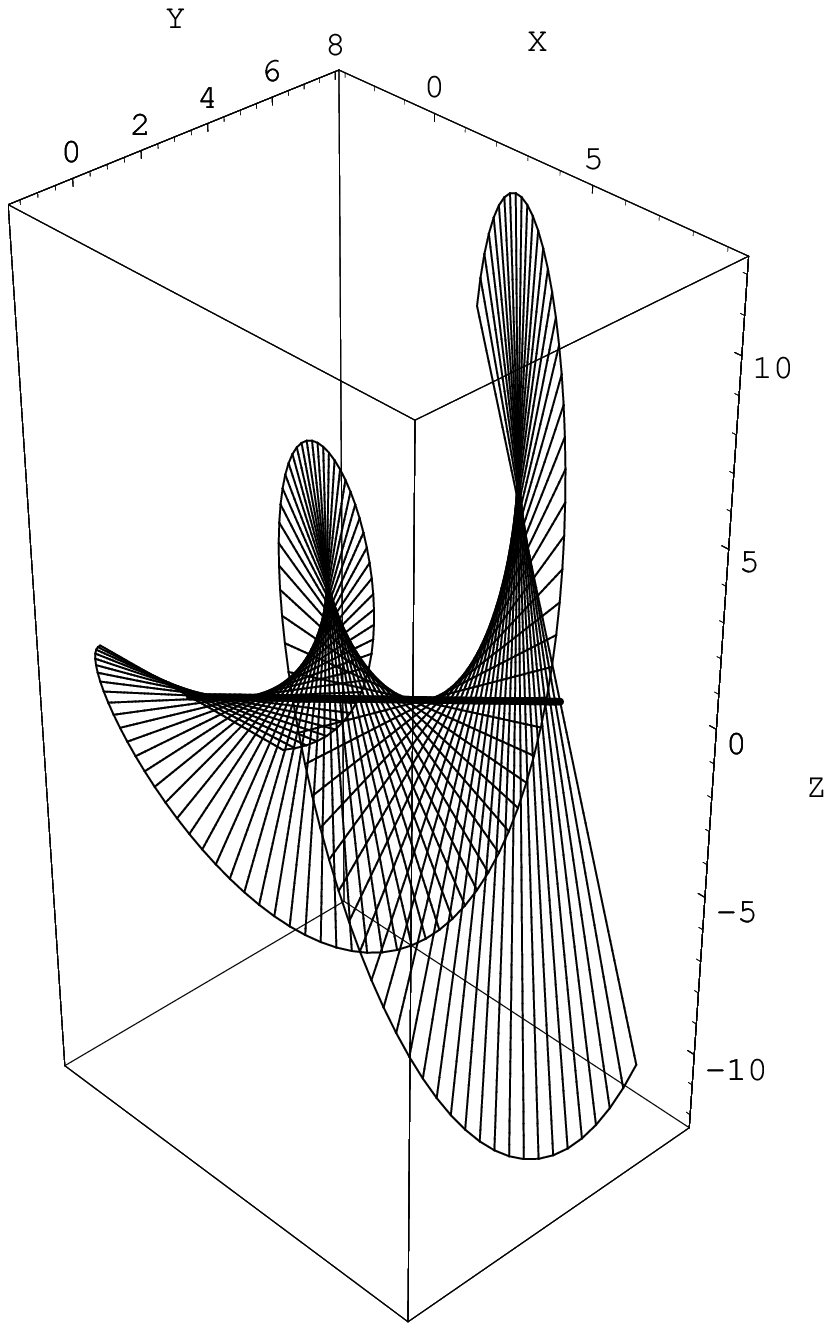}{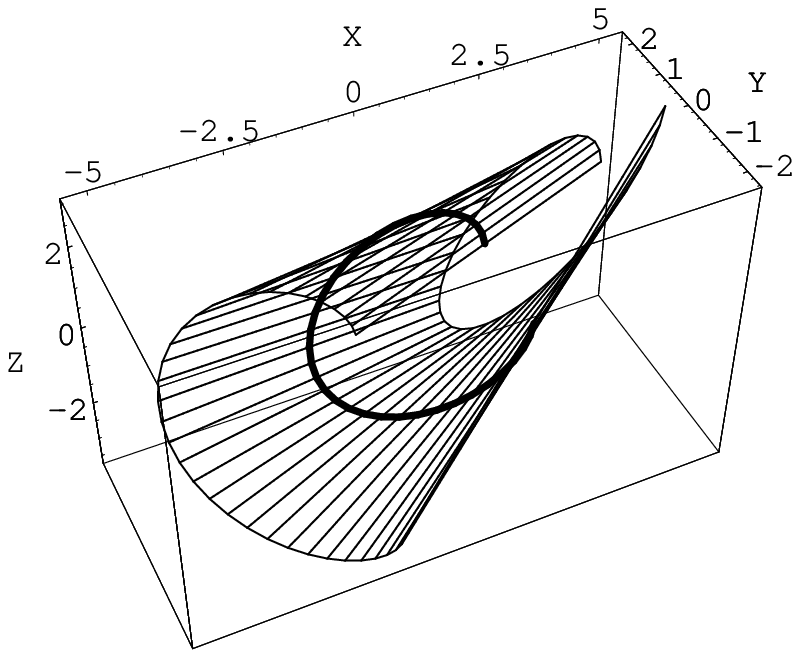}{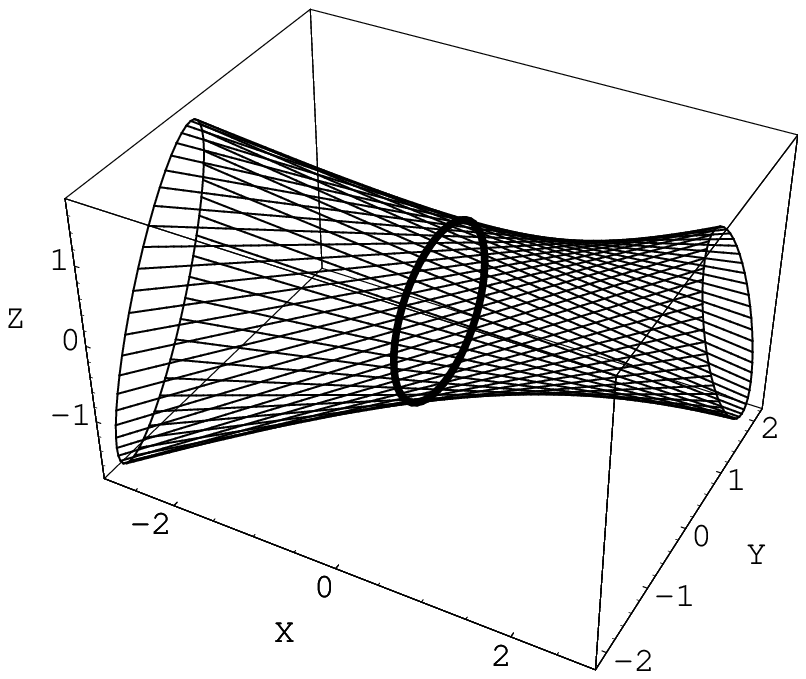}{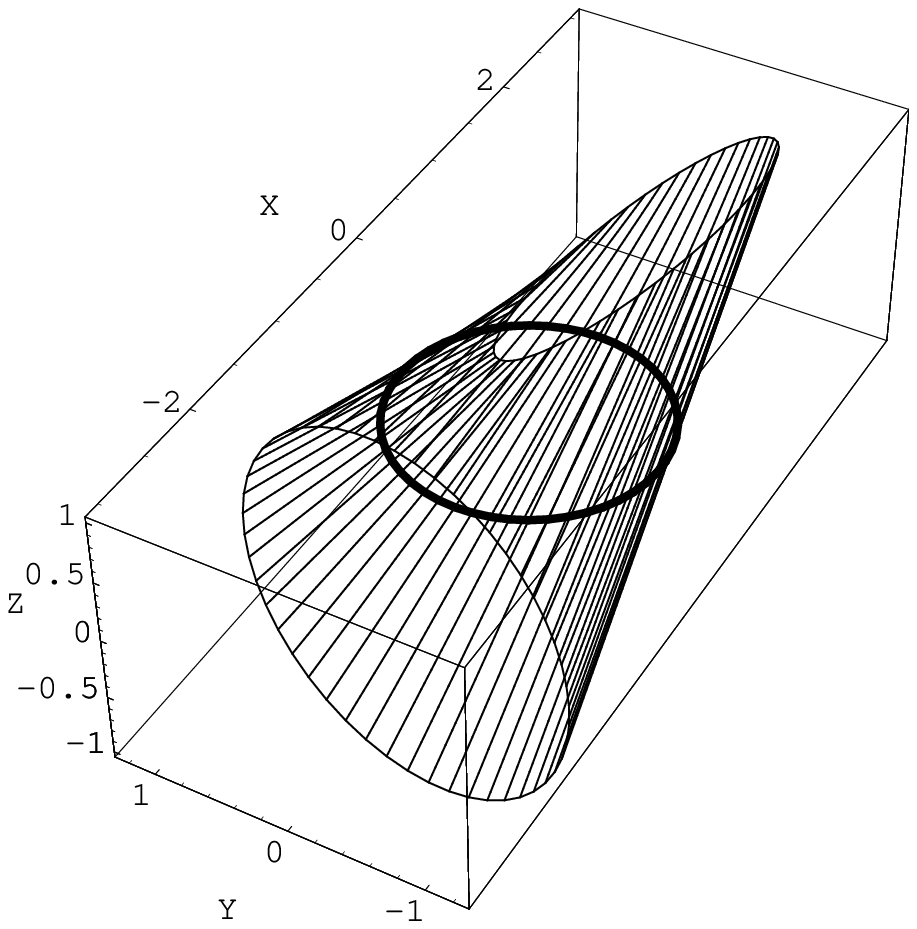}{Examples of $\Delta^{H^1}$-minimal surfaces: a). $\gamma(s)=(0,\cos s,\sin s)$, $\phi(s)=s$, $s\in [0,2\pi]$, 
$t\in[-\frac{3}{2},\frac{3}{2}]$; b) $\gamma(s)=(s,s,0)$, $\phi(s)=s$, $s\in [0,2\pi]$, $t\in [-3,3]$; c). $\gamma(s)=(2\cos s,1-\frac{s}{\pi},2\sin s)$,
$\phi=\frac{\pi}{8}$, $t\in[-1,1]$; d). $\gamma(s)=(0,\cos s,\sin s)$, $\phi=\frac{\pi}{8}$, $t\in[-2,2]$; e). $\gamma(s)=(\cos s,\sin s,0)$, $\phi=\frac{\pi}{25}$, 
$t\in[-2,2]$.}{5.6}{5.6}

\begin{example}{\rm For the Heisenberg-type distribution $\Delta^{H^1}$ condition (\ref{ruledcond}) is always satisfied. 
Thus $\phi$ is constant along characteristic curves.
The characteristic vector field reads 
$$
e^{\phi}=\cos\phi\,\dd_x+\sin\phi\,\dd_y+\frac{1}{2}(x\sin \phi-y\cos\phi)\dd_z,\qquad\phi\in[0,2 \pi].
$$
Thus any characteristic curve satisfies the following system of ODE for some fixed $\phi\in[0,2\pi]$: 
\BE\label{H1sys}
\left\{\begin{array}{lcl}
\dot x&=&\cos\phi\\
\dot y&=&\sin\phi\\
\dot z&=&\frac{1}{2}(x\sin \phi-y\cos\phi)\\
\end{array}\right.
\EE
The solution  of this system is the curve
$$
\left\{\begin{array}{lcl}
x&=&t\,\cos\phi+x_0\\
y&=&t\,\sin\phi+y_0\\
z&=&\frac{1}{2}(x_0\sin \phi-y_0\cos\phi)t+z_0\\
\end{array}\right.,
$$
starting at $(x,y,z)(0)=(x_0,y_0,z_0)$.  We immediately see that the characteristic curves of $\Delta^{H^1}$-minimal surfaces lie on straight lines.
This  fact was first noticed  in \cite{Malc} and it has a lot of important consequences for the global structure of minimal surfaces in $H^1$.
For example, all $\Delta^{H^1}$-minimal surfaces are standard ruled surfaces. The fact that the characteristic curves are straight lines implies that 
any $\Delta^{H^1}$-minimal surface can contain at most one isolated characteristic point.  
In Figure 1 there are shown some examples of  $\Delta^{H^1}$-minimal surfaces, obtained by integration 
of system (\ref{sys}) forward and backward in time with help of Mathematica. The fat line 
denotes the curve of initial conditions $\gamma(s)=q_s(0)$ parametrized by some auxiliary parameter $s$. 
}\end{example}

\begin{example}{\rm In the case of group of roto-translations (see Example 2) condition (\ref{ruledcond}) is satisfied as well, 
and $\phi$ is constant along characteristics. For any fixed $\phi\in[0,2\pi]$ the characteristic vector field is given by
$e^\phi=\cos \phi \cos z \dd_x+\cos\phi\sin z \dd_y+\sin\phi\dd_z$. Let us find explicitly the characteristic curves. They satisfy the following system of ODE:
\BE\label{H1sys}
\left\{\begin{array}{lcl}
\dot x&=&\cos\phi\cos z\\
\dot y&=&\cos\phi\sin z\\
\dot z&=&\sin \phi\\
\end{array}\right.
\EE
The solution that starts at a point $(x_0,y_0,z_0)$  has the form 
$$
\left\{\begin{array}{lcl}
x&=&\frac{\cos\phi}{\sin\phi} \sin z+x_0\\
y&=&-\frac{\cos\phi}{\cos\phi} \sin z+y_0\\
z&=&t \sin\phi+z_0\\
\end{array}\right.\qquad{\rm for}\qquad \phi\ne 0,\pi,
$$
or
$$
\left\{\begin{array}{lcl}
x&=&\pm t\cos z_0+x_0\\
y&=&\pm t\sin z_0+y_0\\
z&=&z_0\\
\end{array}\right.\qquad {\rm for}\qquad \phi= 0,\pi.
$$
In the latter case the characteristic curves are straight lines, though in general this is not true. 
In Figure 2 we present some examples of  $\Delta^{e^2}$-minimal surfaces, constructed by solving numerically equations (\ref{sys}) forward and backward in time by Mathematica.
As before, the fat line denotes the curve of initial conditions $\gamma(s)=q_s(0)$.
\Grarray{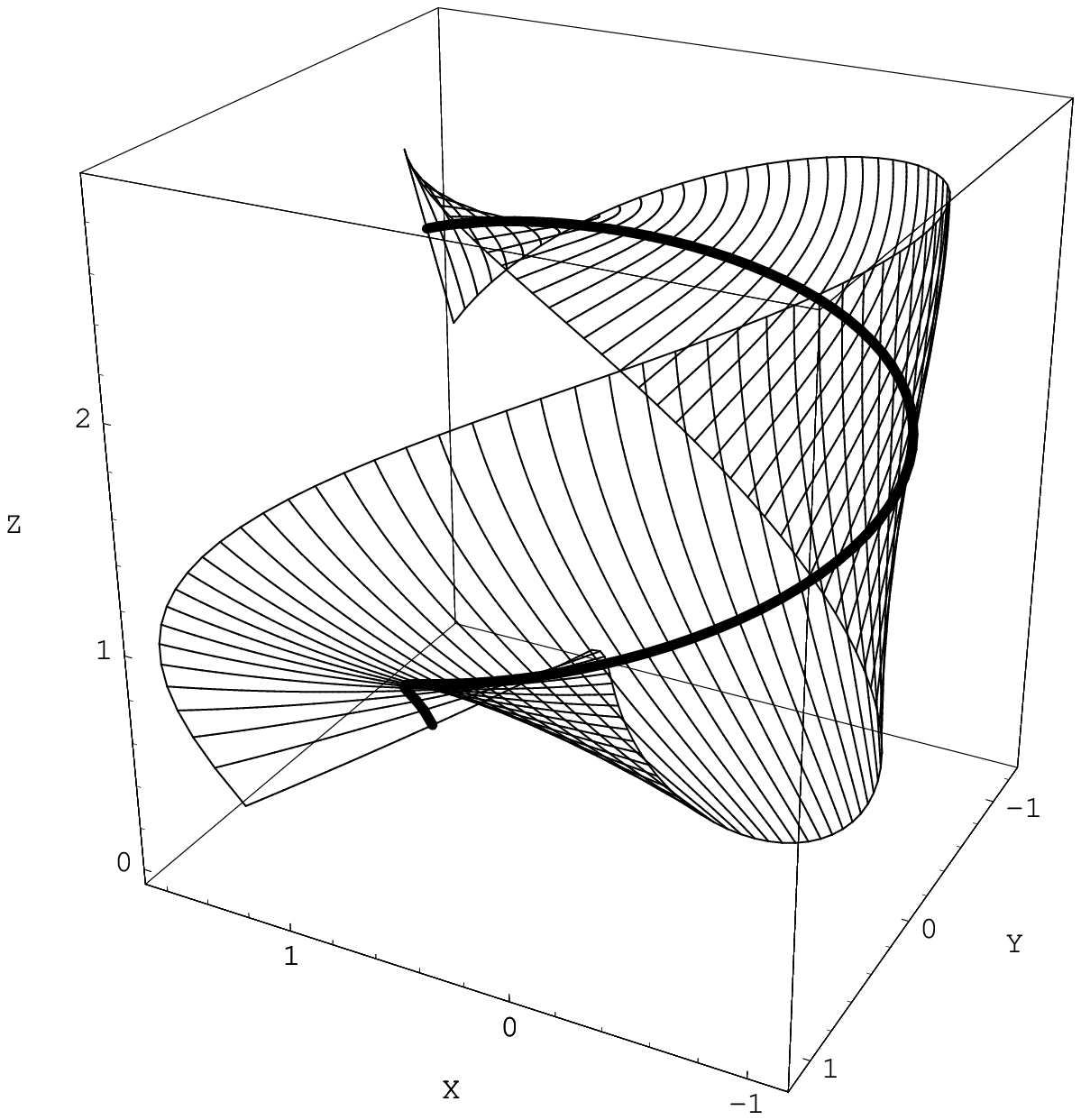}{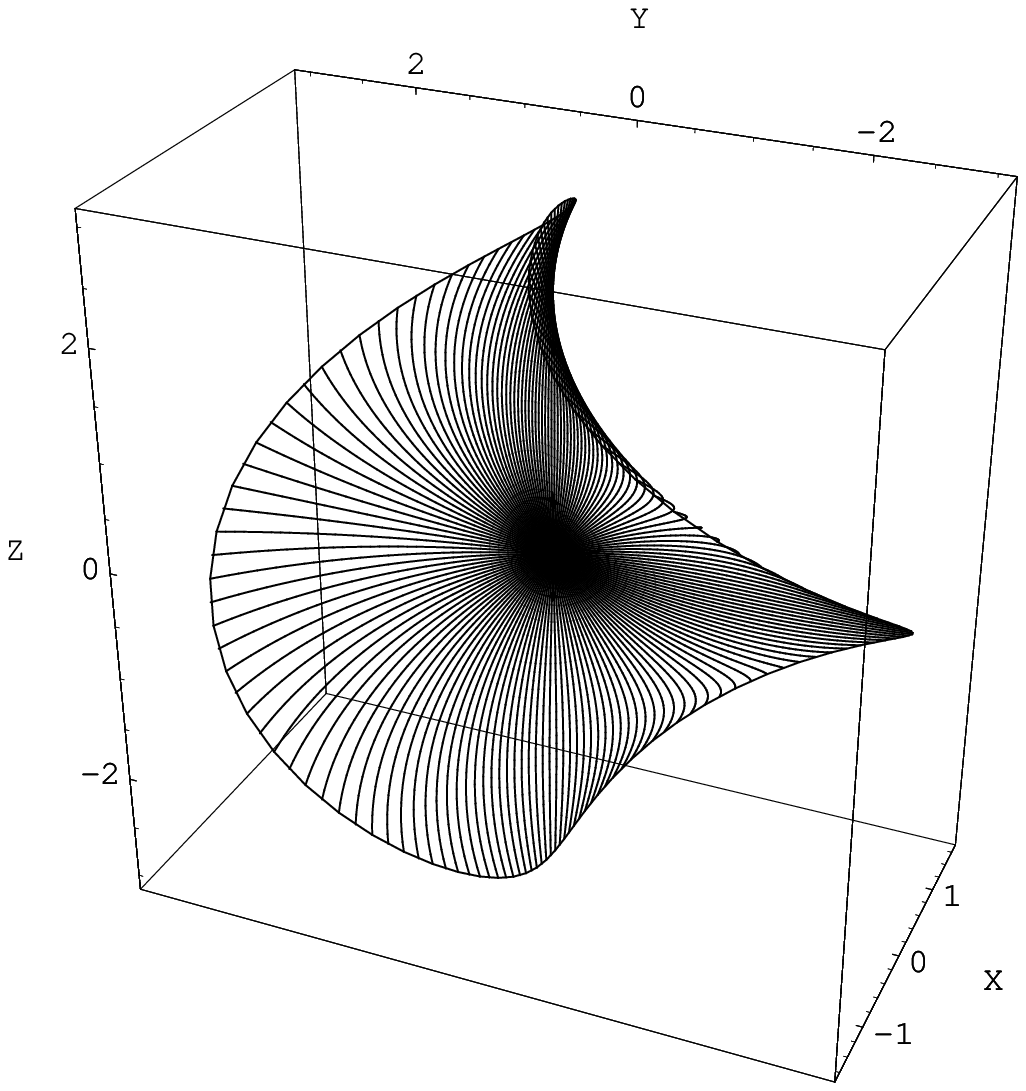}{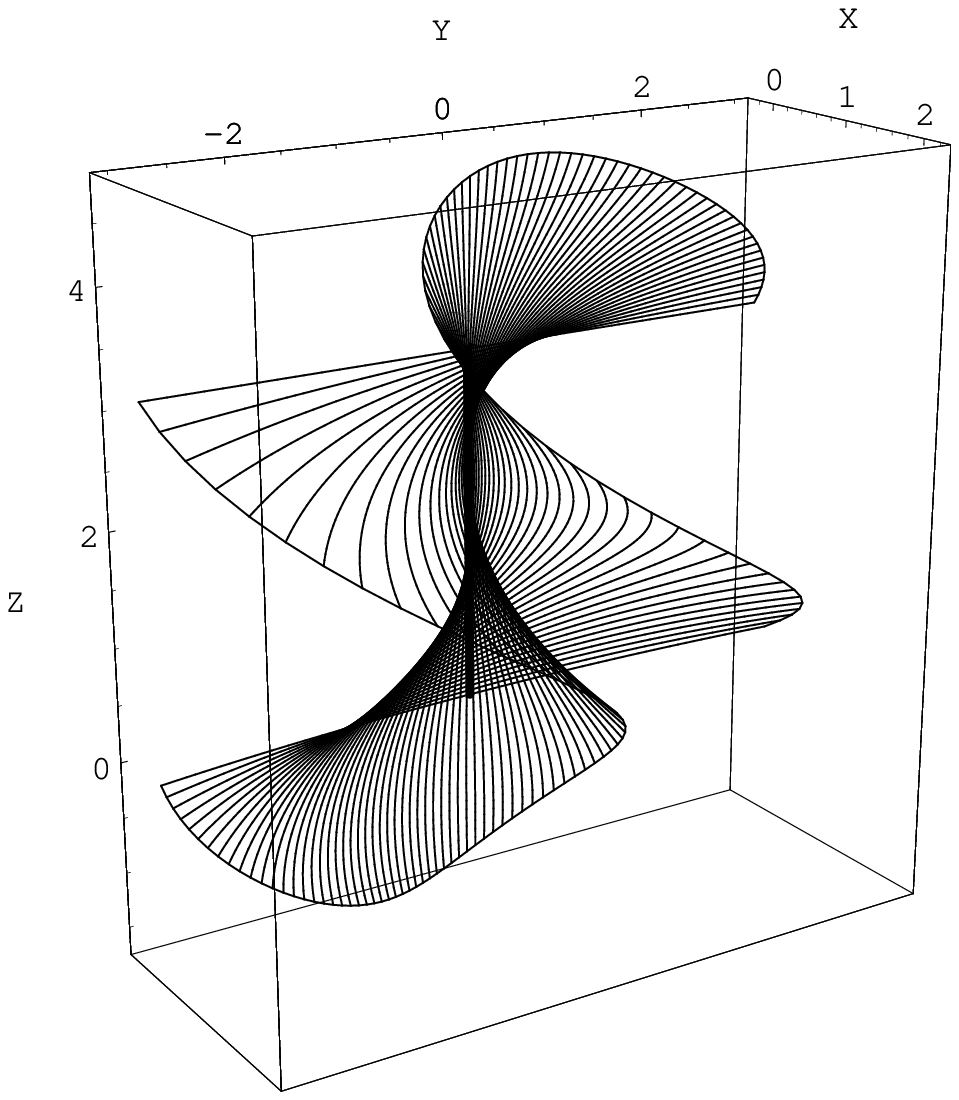}{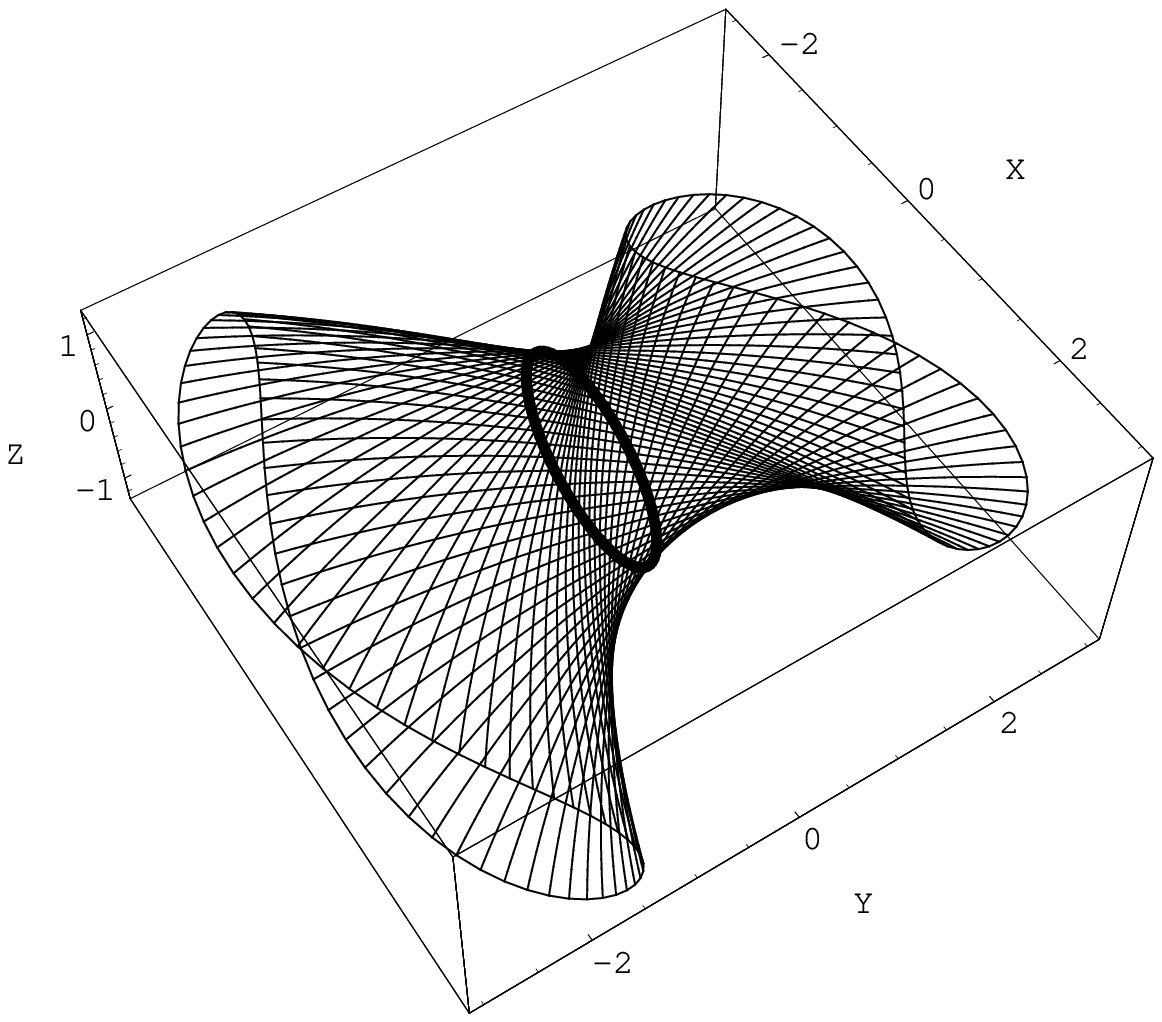}{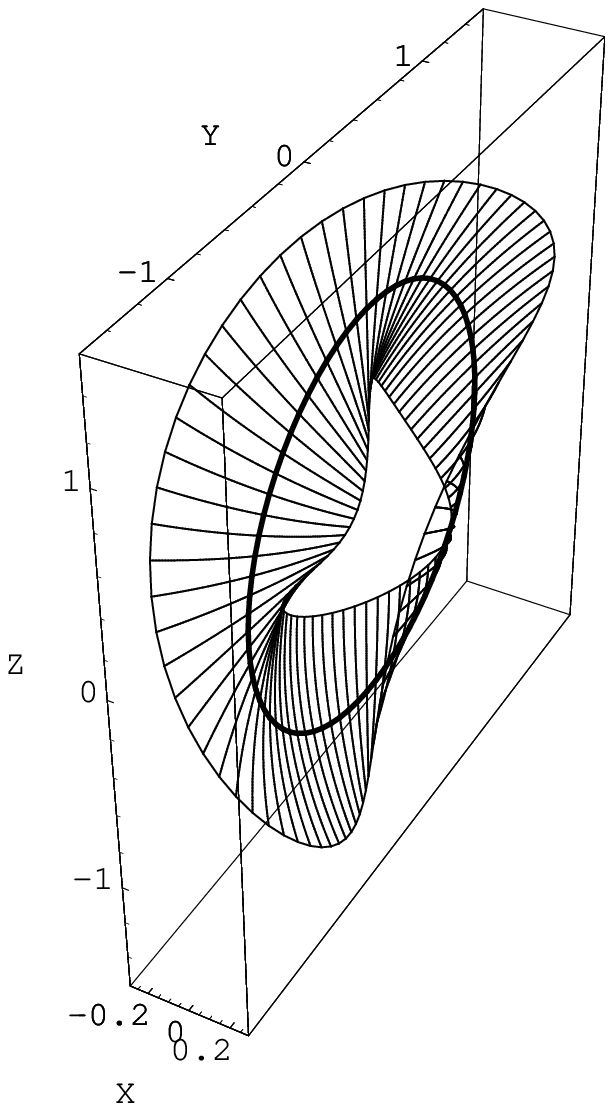}{Examples of $\Delta^{e^2}$-minimal surfaces: a) $\gamma(s)=(\cos s,\sin s,\sqrt s)$, $\phi(s)=\frac{s}{2}$, $s=[0,2\pi]$, 
$t\in[-1,1]$; b) $\gamma(s)=(0,0,0)$, $\phi=\in[0,2\pi]$, $t\in[0,3]$;  c) $\gamma(s)=(0,0,s)$, $\phi(s)=s$, $s\in[0,\pi]$, $t\in[-3,3]$;  d) $\gamma(s)=
(\cos s,0,\sin s)$, $\phi=\frac{2\pi}{75}$, $s\in [0,2\pi]$, $t\in[-3,3]$;  e) $\gamma(s)=(0,\cos s,\sin s)$, $\phi(s)=s$, $s\in [0,2\pi]$, 
$t \in [-\frac{1}{2},\frac{1}{2}]$.}{5.6}{5.6}
}\end{example}

\subsection{Characteristic curves of Sub-Riemannian minimal surfaces and Sub-Riemannian geodesics}
\label{totgeo}

It is natural to ask whether the characteristic curves of Sub-Riemannian minimal surfaces are Sub-Riemannian geodesics.
Recall that Sub-Riemannian geodesics are horizontal curves $t\mapsto \eta(t)$, $\dot\eta(t)\in \Delta_{\eta(t)}$, $t\in [0,T]$, 
which minimize the Sub-Riemannian length
$$
\ell[\eta]=\int\limits_0^T\|\ \dot \eta(\tau)\|_\Delta d\tau
$$
and such that $\|\dot \eta(t)\|_{\Delta}$ is constant for all $t\in[0,T]$ (see \cite{Bel}). The existence of Sub-Riemannian geodesics is guaranteed by the Hopf-Rinow 
theorem provided the distribution  $\Delta$ is bracket-generating. 
According to the classical  Pontryagin Maximum Principle  the Sub-Riemannian geodesics are projections on the base manifold $M$ of the 
corresponding Pontryagin extremals in $T^*M$. In the contact case all these extremals are integral curves of the Hamiltonian vector 
field $\ham \in Vec(T^*M)$\footnote{For the details see \cite{AgrBook}, \cite{AgrExp}.}: 
\BE\label{ham}
\ham=u_1 X_1+u_2 X_2+a\dd_{u_3}-(u_3+b)\dd_{12},
\EE
where
$$
u_i(p,q)=\la p,X_i(q)\ra,\qquad q\in M,\;p\in T^*_q M,\;i=1,2,3,
$$
$$
\dd_{12}=u_1\dd_{u_2}-u_2\dd_{u_1},
$$
$$
a=\ham(u_3)=c_{31}^1(u_1^2-u_2^2)+(c_{32}^1+c_{31}^2)u_1 u_2
$$
$$
b=u_1 c_{12}^1+ u_2 c_{12}^2,
$$
and 
$$
u_1^2+u_2^2=1.
$$
The last condition permits us to introduce a coordinate $\psi\in S^1$ on the oriented circle 
such that  $u_1=\cos\psi$, $u_2=\sin\psi$. Then $\dd_{12}\equiv\dd_{\psi}$.
In particular, it follows that 
\BE\label{hamphi}
\dot\psi=-u_3-b,\qquad\dot u_3=a.
\EE

The characteristic curve of a $\Delta$-minimal surface $W$ is a Sub-Riemannian geodesic if and only if 
\BE\label{pih}
\pi_*(\ham)=e^{\phi},
\EE
i.e., $u_1 X_1+u_2 X_2=e^{\phi}$. Thus we can  identify  $\psi$ and $\phi$. Now comparing (\ref{hamphi}) with (\ref{sys}) we see that the characteristic curves 
are Sub-Riemannian geodesics provided $u_3=0$, which implies $a=0$, i.e., 
\BE\label{cond1}
c_{31}^1\cos 2 \phi+\frac{c_{32}^1+c_{31}^2}{2}\sin 2\phi=0.
\EE
The direct computation yields the solutions of this equation: 
\BE\label{phi*}
\phi^*(q)=\left\{\begin{array}{lc}
-\frac{1}{2}\arctan\left(\frac{2 c_{13}^1(q)}{c_{23}^1(q)+c_{13}^2(q)}\right)+
\frac{k \pi}{2}&{\rm if}\;\;c_{32}^1(q)+c_{31}^2(q)\ne 0;\\
\frac{(2k+1)\pi}{4}& 
\begin{array}{ll}
{\rm if}\;\;c_{32}^1(q)+c_{31}^2(q)=0,\\
\quad\, c_{13}^1(q)\ne 0,\end{array}
\end{array}\right.
\EE
where $k=0,1,2,3$.
The characteristic curves that are Sub-Riemannian geodesics satisfy the ODE $\dot q=e^{\phi^*(q)}(q)$.
If both coefficients in (\ref{cond1}) are zero, then all characteristic curves are Sub-Riemannian geodesics, parametrized according to (\ref{hamphi}).

\vspace{10pt}

In the particular case of  Lie groups, the structural constants $c_{ij}^k$  do not depend on the point of the base manifold. Hence $\phi^*=const$ 
and we have an additional condition:
\BE\label{cond2}
b=c_{12}^1\cos\phi+ c_{12}^2\sin\phi=0.
\EE
In general, this condition is stronger than (\ref{cond1}).
It is easy to check that if (\ref{cond2}) is  non-degenerate, its solutions belong to the set of solutions of (\ref{cond1}).
Indeed, if $c_{12}^2\ne 0$, then combining (\ref{cond1}) and (\ref{cond2})
by direct computation we obtain the following compatibility condition for the structural constants:
\BE\label{sconst}
c_{13}^1((c_{12}^1)^2-(c_{12}^2)^2)+c_{12}^1c_{12}^2(c_{23}^1+c_{13}^2)=0.
\EE
Comparing (\ref{sconst}) with (\ref{jac}) one can see that it is trivially satisfied.
If $c_{12}^1\ne 0$, $c_{12}^2=0$, then  (\ref{jac}) implies $c_{13}^1=c_{13}^2=0$ and (\ref{cond1}) reduces to $c_{23}^1 \sin 2\phi=0$.
On the other hand, in this case (\ref{cond2}) yields $\cos \phi=0$, which clearly satisfy (\ref{cond1}). 
Summing up, we obtain the following classification:\\
a). if $c_{12}^2\ne 0$, then $\phi^*=-\arctan\left(\frac{c_{12}^1}{ c_{12}^2}\right)+k\pi$, $k=0,1$;\\
b). if $c_{12}^2=0$, then $\phi^*=\frac{k\pi}{2}$, $k=1,3$;\\
c). if $c_{12}^1=c_{12}^2=0$  then $\phi^*$ is given by (\ref{phi*});\\
d). all characteristic curves are Sub-Riemannian geodesics if both (\ref{cond1}) and (\ref{cond2}) degenerate. \\
We conclude this discussion by analysis of the characteristic curves in $H^1$ and $e^2$. 

\begin{example}{\rm In the Heisenberg case, since all structural constants, but $c_{12}^3$, vanish, the parameter $\phi^*$ can take any value in 
$[0,2\pi]$. As we have already seen, the characteristic curves are straight lines and they all are Sub-Riemannian geodesics.
Fixing a point $q_0\in M$ and varying $\phi^*\in [0,2\pi]$ one can  generate a plane, which is a totally geodesic surface in the Sub-Riemannian sense, 
and in the same time it is an entire  $\Delta^{H^1}$-minimal surface with one characteristic point  at $q_0$.
}\end {example}

\begin{example}{\rm In the case of the  distribution $\Delta^{e^2}$ 
there are two non-zero structural constants $c_{12}^3=c_{23}^1=1$.
From (\ref{cond1}) we obtain $\phi^*=\frac{k\pi}{2}$, $k=1,3$.  
A simple calculation shows that there are two families of characteristic curves that are Sub-Riemannian geodesics:
$$
x=t\,\cos z_0+x_0,\qquad y=t\,\sin z_0+y_0\qquad z=z_0;
$$
$$
x=x_0,\qquad y=y_0,\qquad z=\pm t+z_0.
$$
Note that these curves are actually  straight lines.
}
\end {example}

\end{document}

%% file: newcom1.tex
\newcommand{\real}{\mathbb{R}}

\newcommand{\eps}{\varepsilon}

\newcommand{\intt}{\int\limits}
\newcommand{\ham}{\vec h} 
 
\newcommand{\dd}{\partial}

\newcommand{\Grarray}[8]
{
\begin{figure}\begin{center}
a)~\includegraphics[width=#7truecm,height=#8truecm]{#1.eps}
b)~\includegraphics[width=#7truecm,height=#8truecm]{#2.eps}
c)~\includegraphics[width=#7truecm,height=#8truecm]{#3.eps}
d)~\includegraphics[width=#7truecm,height=#8truecm]{#4.eps}
e)~\includegraphics[width=#7truecm,height=#8truecm]{#5.eps}
\caption{#6}
\label{#1}
\end{center}
\end{figure}
\noindent$\!$}

\newtheorem{theorem}{Theorem}

\newtheorem{defi}{Definition}

\newcommand{\BE}{\begin {equation}}
\newcommand{\EE}{\end {equation}}